\documentclass[10pt,reqno,final]{article}

\usepackage{amsmath,amsfonts,amssymb,amsthm,version}
\usepackage{mathrsfs,fancybox,pifont}
\usepackage{graphicx}
\usepackage{url,hyperref}
\usepackage[notcite,notref]{showkeys}
\usepackage{color}
\usepackage{subfigure,multirow}
\usepackage{epstopdf}
\usepackage{cases}
\usepackage{mathtools}
\usepackage{algorithm,algorithmic}
\usepackage{authblk}
\usepackage{lipsum}
\usepackage{cite}

\allowdisplaybreaks

\theoremstyle{plain}
\numberwithin{equation}{section}
\newtheorem{lemma}{Lemma}[section]
\newtheorem{theorem}{Theorem}[section]

\newtheorem{remark}{Remark}[section]

\newtheorem{problem}{Problem}
\newtheorem{assumption}{Assumption}[section]

\theoremstyle{definition}

\title{Error analysis of an effective numerical scheme for a temporal multiscale plaque growth problem}
\author{Xinyu Li\thanks{School of Mathematics and Physics, University of Science and Technology Beijing, Beijing 100083, China} \,,
Zhaoyang Wang\thanks{Corresponding author. Department of Applied Mathematics, University of Science and Technology Beijing, Beijing 100083, China (zhaoyang584520@163.com)} \,,
Mingjie Liao\thanks{AsiaInfo Technologies (China) Inc, Beijing, China (liaomj5@asiainfo.com)} \, and
Ping Lin\thanks{Corresponding author. Division of Mathematics, University of Dundee, Dundee DD1 4HN, United Kingdom (p.lin@dundee.ac.uk)}}

\affil{}

\date{}

\begin{document}
\maketitle

\begin{abstract}
In this work, we propose a simple numerical scheme based on a fast front-tracking approach for solving a fluid-structure interaction (FSI) problem of plaque growth in blood vessels. A rigorous error analysis is carried out for the temporal semi-discrete scheme to show that it is first-order accurate for all macro time step $\Delta T$, micro time step $\Delta t$ and scale parameter $\epsilon$. A numerical example is presented to verify the theoretical results and demonstrate the excellent performance of the proposed multiscale algorithm.


\medskip
\end{abstract}

\section{Introduction}
\label{section1}
Many problems in biology have multiscale characteristics. For example, in the formation of atherosclerotic plaque in arteries, blood flow on the scale of milliseconds to seconds has an effect on the plaque growth in the vessel, while plaque formation takes a long time, from months to several years. 

\begin{figure}[H]
\centering
\includegraphics[width=9cm,height=2.5cm]{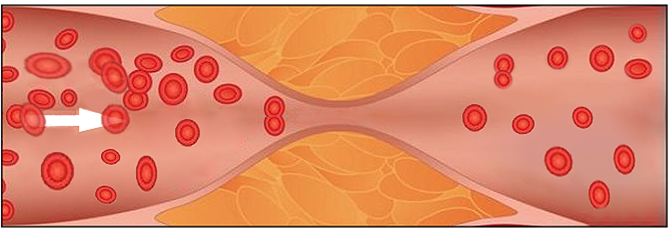}
\caption{Schematic diagram of plaque growth causing narrowing of blood vessels.}
\label{plaque}
\end{figure}

In this work, we are interested in the long-term simulation of the following fluid-structure interaction of plaque growth (See \cite{frei2020efficient} and a schematic diagram in Figure \ref{plaque}):
\begin{problem}
\label{problem1}
\begin{equation}\label{the1.1}
\begin{split}
&\nabla\cdot \textbf{v}=0, \quad \rho(\frac{\partial\textbf{v}}{\partial t}+(\textbf{v}\cdot\nabla)\textbf{v})=\text{div} \ \sigma(\textbf{v},p)+\textbf{f}, \quad \text{in} \ \Omega(u(t)), \ \textbf{v}(0)=\textbf{v}_0,\\
&\frac{d}{dt}u(t)=\epsilon R(\textbf{v},u), \ u(0)=u_0,\\
&\Omega(u(t))=\left\{(x,y):|x|<a, |y|<b-\gamma(u,x)\right\}.
\end{split}
\end{equation}
\end{problem}
The velocity $\textbf{v}$ and concentration $u$ represent the micro (fast) variable and the macro (slow) variable, respectively. $\rho$ is the density of blood which we assume to be constant in this paper, $\sigma=-pI+\rho\nu(\nabla\textbf{v}+\nabla\textbf{v}^T)$ the Cauchy stress with the kinematic viscosity $\nu$. A periodic force $\textbf{f}(t)=\textbf{f}(t+1)$ is applied to the flow due to the periodic nature of the heart pulse. $\epsilon\ll 1$ is a small scale parameter that controls the change of $u$. The reaction term $R\geq 0$ describes the influence of wall shear stress on the boundary growth, which can be expressed as follows \cite{frei2020efficient}
\begin{equation}\label{the1.2}
\begin{split}
&R=(1+u)^{-1}(1+|\sigma_{WSS}(\textbf{v})|^2)^{-1}=O(1),\ \sigma_{WSS}(\textbf{v})=\sigma_0^{-1}\int_{\varGamma_{\text{wall}}}\rho\nu(I-\vec{n}\vec{n}^{T})(\nabla\textbf{v}+\nabla\textbf{v}^T)\vec{n}ds,
\end{split}
\end{equation}
where $\vec{n}$ denotes the outward facing unit normal vector at the deformation boundary $\partial\Omega$.

For this plaque growth problem with temporal multiscale features, the micro scale has to be resolved by time steps in the order of centiseconds \cite{figueroa2009computational}. Hence, it is computationally expensive to do the realistic numerical simulation of plaque growth over several months.  

An effective approximation for problem \ref{problem1} was introduced in \cite{frei2020efficient} based on the arbitrary Lagrangian-Eulerian (ALE) framework, but the error between the approximation problem and the original problem is only analyzed for a largely simplified blood flow/plaque growth (ODE/ODE) system. In a recent work \cite{wang2022fast} an effective approximation for a more complicated plaque growth problem with memory effects (fractional growth) was presented based on a simple but fast front-tracking framework, and the approximation analysis was done for the original flow/plaque growth (PDE/fractional ODE) system. However, a full error analysis of the discrete approximation problem, especially for the periodic flow problem involved in the effective approximation was not done in both above mentioned papers. In addition, to our knowledge, the numerical error analysis of a time-peridic flow model has not been considered much in literature.

In this paper, our goal is to carry out a complete error analysis for a temporal semi-discrete scheme of this multiscale plaque growth problem. In addition, the numerical implementation is based on the open-source finite element Software FreeFEM++ \cite{hecht2012new}, which demonstrates its applicability for evaluating multiphysics models involving small deformations.

The rest of paper is organized as follows. In Section \ref{section2}, we introduce an effective auxiliary problem that approximates problem \ref{problem1}, and based on this, we design temporal numerical schemes for both micro and macro parts of the approximate problem and formulate an effective multiscale algorithm. In Section \ref{section3}, we provide rigorous error estimate for this temporal semi-discrete scheme. Numerical results are presented in Section \ref{section4}.

\textbf{Notation.} For space $\Omega$ and $m\geq 0$, we use the standard notation for the Sobolev space $H^{m}(\Omega)$ and the Banach space $L^m(\Omega)$. We denote spaces $H^{m}(\Omega)$ and $L^{m}(\Omega)$ by $H^{m}$ and $L^{m}$ in short. We denote by $(\cdot,\cdot)$ the inner product, and by $C$ any constant independent of $\epsilon$ and $t$.

\section{Multiscale algorithm and implementation}
In this section, we present a time-periodic auxiliary problem, design an 
effective multiscale algorithm based on it, and provide implementation details.
\label{section2}
\subsection{A time-periodic auxiliary problem}
In order to use the fast multiscale algorithm, the following time-periodic auxiliary problem with homogeneous Dirichle boundary condition is introduced:
\begin{problem}
\label{problem2}
\begin{subequations}\label{the2.1}
\begin{align}
&\nabla\cdot \textbf{v}_U=0,
\quad \rho\left(\frac{\partial\textbf{v}_U}{\partial t}+(\textbf{v}_U\cdot\nabla)\textbf{v}_U\right)=\nu \Delta \textbf{v}_U-\nabla p_U+\textbf{f} \quad \text{in} \ [0,1]\times \Omega(U), \label{the2.1a}\\
&\textbf{v}_U(0)=\textbf{v}_U(1) \quad \text{in} \ \Omega(U),\\
&\frac{d}{dt}U(t)=\epsilon\int_{0}^{1} R\left(\textbf{v}_{U(t)}(s), U(t)\right)ds, \quad U(0)=u_0, \label{the2.1c}
\end{align}
\end{subequations}
\end{problem}
Here, $U(t)=\int_{t}^{t+1}u(s)ds$ was defined in \cite{frei2020efficient}.

By following arguments and techniques in \cite{frei2020efficient} and \cite{wang2022fast} it can be proved that problem \ref{problem2} can be regarded as an approximate problem of problem \ref{problem1} and the approximation error $|U(t)-u(t)|=O(\epsilon)$. Namely, the solution of the auxiliary problem \ref{problem2} can be used as an approximation of the solution of the original problem \ref{problem1}. 

Different from the original problem \ref{problem1}, the macro variable $U$ of the auxiliary problem \ref{problem2} is homogenized by time and the micro equation uses time-periodic condition instead of the initial value condition. This makes it possible to use a large time step to approximate the macro variable $U$. At each of these macro-time-steps, the micro equation is initialized based on the current macro variable and can be solved on the interval $[0,1]$ due to its periodicity 1. 

\begin{remark}
The error analysis of approximation problem \ref{problem2} given in \cite{wang2022fast} is based on the linearized Stokes problem. The analysis for full Navier-Stokes equations remains an interesting open problem.
\end{remark}

\subsection{Numerical schemes and algorithms}
Since $R=O(1)$, it is not difficult to show $|u(t)|=O(\epsilon t)$. So to observe an $O(1)$ growth, we need to compute up to a time $T=O(\epsilon^{-1})$. We first split the time interval $[0,T=O(\epsilon^{-1})]$ into subintervals of equal size $0=T_0<T_1<\cdots<T_m<\cdots<T_M=T$, and the macro step size $\Delta T=T_m-T_{m-1}$. Then we split the fast periodic interval $[0,1]$ into subintervals of equal size $0=t_0<t_1<\cdots<t_n<\cdots<t_N=1$, and the micro step size $\Delta t=t_n-t_{n-1}$. 

Based on (\ref{the2.1a})-(\ref{the2.1c}), we propose an explicit discrete scheme of U:
\begin{equation}\label{the2.2}
\begin{split}
U_{m}=U_{m-1}+\Delta T\epsilon\sum\limits_{n=1}^{N}\Delta t R\left((\textbf{v}_{U_{m-1}})_n, U_{m-1}\right), \quad  U_0=u_0.
\end{split}
\end{equation}

We approximate the time-periodic solution $\textbf{v}_U$ on the micro scale with the simple Euler scheme:
\begin{equation}\label{the2.3}
\begin{split}
&\nabla\cdot (\textbf{v}_{U_m})_{n}=0, \quad \rho\left(\frac{(\textbf{v}_{U_m})_{n}-(\textbf{v}_{U_m})_{n-1}}{\Delta t}+((\textbf{v}_{U_m})_{n}\cdot\nabla)(\textbf{v}_{U_m})_{n}\right)=\nu \Delta (\textbf{v}_{U_m})_{n}-\nabla (p_{U_m})_{n}+\textbf{f}_n,\\
&\textbf{v}_{U_m}(0)=\textbf{v}_{U_m}(1), \quad \text{in} \ \Omega(U_m).
\end{split}
\end{equation}

It can be seen that for a fixed $U$ (flow domain is fixed), we need to solve the time-periodic  Navier-Stokes equations without an initial value. Two different iterative algorithms are designed in \cite{frei2020efficient} and \cite{wang2022fast} to solve for the time-periodic flow equations. Here we apply the iterative algorithm developed in \cite{wang2022fast} (See Algorithm 4.2) for identifying the initial value of the periodic-in-time solution to the Navier-Stokes equations. 

\begin{figure}[H]
	\centering
	\includegraphics[width=14cm,height=5cm]{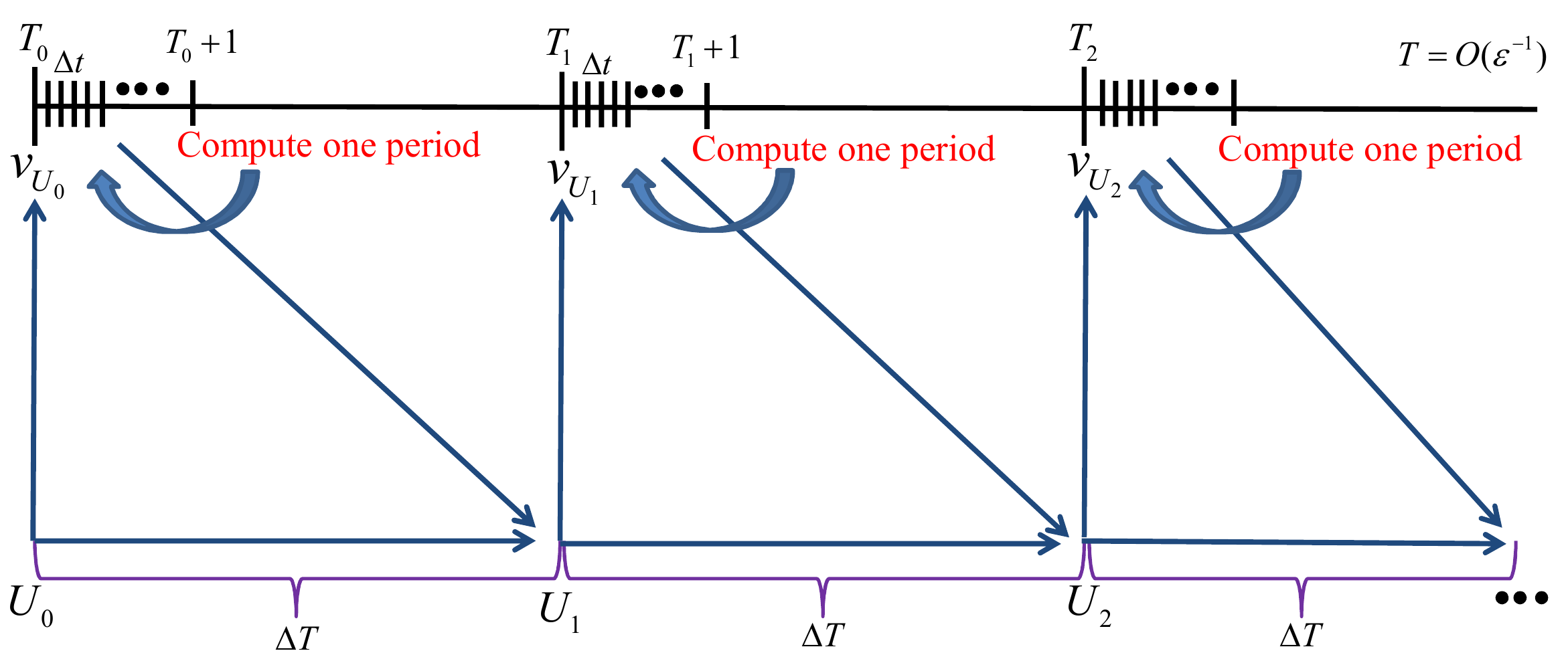}
	\caption{Schematic diagram of the multiscale method. Stepping the macro variable $U(t)$ with macro step size $\Delta T\gg 1$, then computing the time-periodic solution $\textbf{v}_{U}(t)$ on one cycle by using the micro step size $\Delta t$ for a fixed flow domain $\Omega(U)$.}
	\label{Multiscale}
\end{figure}

Before presenting our multiscale algorithm, we first illustrate the effective numerical framework in Figure \ref{Multiscale}. We would like to remark that our proposed multiscale algorithm is based on a simple front-tracking method, which uses the explicit scheme (\ref{the2.2}) to determine the computational domain and then solves equation (\ref{the2.3}). This is significantly different from the ALE approach used in \cite{frei2020efficient} which results in a more complicated ALE transformed Navier-Stokes equations.

Now, we write down the following effective numerical multiscale algorithm.

\begin{algorithm}[H]
\caption{Multiscale front-tracking algorithm.}
\label{algorithm1}
\begin{algorithmic}[1]
\STATE{\textbf{Input:} Set $U_0=u_0$. Given an inflow velocity $\textbf{v}_0$, let $0<\tau \ll 1$ be a given tolerance for identifying time-periodic solutions.}
\STATE{\textbf{Output:} $U_1, U_2,..., U_M$.}
\FOR{$m=1$ to $M$} 
\STATE{Compute equation (\ref{the2.3}) by the iterative algorithm 2 in \cite{wang2022fast} to obtain $(\textbf{v}_{U_{m-1}})_0$, $(\textbf{v}_{U_{m-1}})_1$, $(\textbf{v}_{U_{m-1}})_2$,...,$(\textbf{v}_{U_{m-1}})_n$}
\STATE{Compute the reaction term $R_{m-1}$ in equation (\ref{the2.2})
$$R_{m-1}=\Delta t\sum\limits_{n=1}^{N} R\left((\textbf{v}_{U_{m-1}})_n, U_{m-1}\right)$$}
\STATE{Step forward $U_{m-1}\to U_{m}$ with the explicit scheme
$$U_m=U_{m-1}+\Delta T\epsilon R_{m-1}$$}
\ENDFOR
\end{algorithmic}
\end{algorithm}

\section{Error analysis of the effective temporal semi-discrete scheme}
\label{section3}
In this section, we analyse the error of the temporal semi-discrete scheme for the macro equation based on (\ref{the2.2}) and (\ref{the2.3}).

We first consider the error estimate of the micro scale.

\begin{assumption}
\label{assumption}
Let $U$ be fixed and $\Omega(U)$ be a domain with a boundary that is sufficiently smooth or polygonal with a finite number of convex corner. We assume that the solution $(\textbf{v}_{U},p_U)$ of (\ref{the2.1}) is unique, and satisfies the following regularities:
\begin{equation}\label{the3.1}
\begin{split}
\textbf{v}_U\in L^{\infty}(0,T;H^2(\Omega(U))), \ d_t\textbf{v}_U\in L^{2}(0,T;L^2(\Omega(U))), \ p_U\in L^{\infty}(0,T;H^1(\Omega(U))).
\end{split}
\end{equation}
\end{assumption}

\begin{lemma}
\label{lemma}
Under the assumption \ref{assumption}, for a fixed flow domain $\Omega(U_m)$, the error of the micro scale scheme (\ref{the2.3}) can be expressed as follows
\begin{equation}\label{the3.2}
\begin{split}
\Vert \textbf{v}_{U_m}(t_n)-(\textbf{v}_{U_m})_n\Vert_{H^1}^2+\Delta t\sum\limits_{n=1}^{N}\Vert \textbf{v}_{U_m}(t_n)-(\textbf{v}_{U_m})_n\Vert_{H^2}^2\leq C_{L31}(\Delta t)^2
\end{split}
\end{equation}
for all $1\leq n\leq N$. Here, for simplicity of analysis we consider homogeneous Dirichlet boundary condition.
\begin{proof}
For a fixed $U_m$ and $1\leq n\leq N$, we take $t=t_{n}$ in the Navier-Stokes equations of (\ref{the2.1}) to obtain that
\begin{equation}\label{the3.3}
\begin{split}
\rho\left(\frac{\textbf{v}_{U_m}(t_n)-\textbf{v}_{U_m}(t_{n-1})}{\Delta t}+(\textbf{v}_{U_m}(t_n)\cdot\nabla)\textbf{v}_{U_m}(t_n)\right)=\nu \Delta \textbf{v}_{U_m}(t_n)-\nabla p_{U_m}(t_n)+\textbf{f}(t_n)+r_n
\end{split}
\end{equation}
with $\nabla\cdot \textbf{v}_{U_m}(t_n)=0$. Under the regularity assumption \ref{assumption}, we have $\Delta t\sum\limits_{n=1}^{N}\Vert r_n\Vert_{L^2}^2\leq C(\Delta t)^2$.

Let $\textbf{e}_n=\textbf{v}_{U_m}(t_n)-(\textbf{v}_{U_m})_n$. The error equation is 
\begin{equation}\label{the3.4}
\begin{split}
\rho\left(\frac{\textbf{e}_n-\textbf{e}_{n-1}}{\Delta t}\right)+((\textbf{v}_{U_m})_n\cdot \nabla)\textbf{e}_n+(\textbf{e}_n\cdot \nabla)\textbf{v}_{U_m}(t_n)=\nu \Delta \textbf{e}_n-\nabla p_{U_m}(t_n)+\nabla (p_{U_m})_{n}+r_n
\end{split}
\end{equation}
with $\nabla\cdot \textbf{e}_n=0$.

Multiplying (\ref{the3.4}) by $2\Delta t\textbf{e}_n$ and integrating over $\Omega(U_m)$, one has
\begin{equation}\label{the3.5}
\begin{split}
\rho\left(\Vert \textbf{e}_n\Vert_{L^2}^2-\Vert \textbf{e}_{n-1}\Vert_{L^2}^2+\Vert \textbf{e}_n-\textbf{e}_{n-1}\Vert_{L^2}^2\right)+2\nu\Delta t\Vert \nabla\textbf{e}_n\Vert_{L^2}^2=-2\Delta t\left((\textbf{e}_n\cdot \nabla)\textbf{v}_{U_m}(t_n),\textbf{e}_n\right)+2\Delta t(r_n,\textbf{e}_n).
\end{split}
\end{equation}

By using the H\"older inequality, the Soblev inequality and the Young inequality, we estimate the right-hand side of (\ref{the3.5})
\begin{equation}\label{the3.6}
\begin{split}
\text{RHS}&\leq C\Delta t\Vert \textbf{e}_n\Vert_{L^2}\Vert \nabla\textbf{v}_{U_m}(t_n)\Vert_{L^\infty}\Vert \textbf{e}_n\Vert_{L^{2}}+\Delta t\Vert r_n\Vert_{L^2}^2+\Delta t\Vert \textbf{e}_n\Vert_{L^2}^2\\
&\leq C\Delta t\Vert \textbf{e}_n\Vert_{L^2}^2+\Delta t\Vert r_n\Vert_{L^2}^2.
\end{split}
\end{equation}

Summing up (\ref{the3.5}) from $1$ to $n$, and using the discrete Gronwall inequality, we obtain
\begin{equation}\label{the3.7}
\begin{split}
\Vert \textbf{e}_k\Vert_{L^2}^2+\sum\limits_{k=1}^{n}\Vert \textbf{e}_k-\textbf{e}_{k-1}\Vert_{L^2}^2+\Delta t\sum\limits_{k=1}^{n}\Vert \nabla\textbf{e}_k\Vert_{L^2}^2\leq C(\Delta t)^2.
\end{split}
\end{equation}

Next, we estimate $(\textbf{v}_{U_m})_n$ in $H^2$-norm. For $\Vert (\textbf{v}_{U_m})_n-(\textbf{v}_{U_m})_{n-1}\Vert_{L^2}$, we have
\begin{equation}\label{the3.8}
\begin{split}
\Vert (\textbf{v}_{U_m})_n-(\textbf{v}_{U_m})_{n-1}\Vert_{L^2}^2&\leq 2\Vert \textbf{e}_n-\textbf{e}_{n-1}\Vert_{L^2}^2+2\Vert_{L^2}^2+2\Vert \textbf{v}_{U_m}(t_n)-\textbf{v}_{U_m}(t_{n-1})\Vert_{L^2}^2\\
&\leq 2\Vert \textbf{e}_n-\textbf{e}_{n-1}\Vert_{L^2}^2+2\Delta t\int_{t_{n-1}}^{t_{n}}\Vert d_t\textbf{v}_{U_m}(t)\Vert_{L^2}^2dt.
\end{split}
\end{equation}
It follows that 
\begin{equation}\label{the3.9}
\begin{split}
\Delta t\sum\limits_{n=1}^{N}\left\Vert\frac{ (\textbf{v}_{U_m})_n-(\textbf{v}_{U_m})_{n-1}}{\Delta t}\right\Vert_{L^2}^2\leq C.
\end{split}
\end{equation}

By using the regularity of the solution to the Stokes problem, we have
\begin{equation}\label{the3.10}
\begin{split}
\Delta t\sum\limits_{n=1}^{N}(\Vert(\textbf{v}_{U_m})_n\Vert_{H^2}^2+\Vert (p_{U_m})_n\Vert_{H^1}^2)\leq \Vert\textbf{f}\Vert_{L^2}^2+\Delta t\sum\limits_{n=1}^{N}\left\Vert\frac{ (\textbf{v}_{U_m})_n-(\textbf{v}_{U_m})_{n-1}}{\Delta t}\right\Vert_{L^2}^2+\Vert \nabla(\textbf{v}_{U_m})_n\Vert_{L^2}^2\leq C.
\end{split}
\end{equation}

Multiplying (\ref{the3.4}) by $2\Delta t(\textbf{e}_n-\textbf{e}_{n-1})$ and integrating over $\Omega(U_m)$, we obtain
\begin{equation}\label{the3.11}
\begin{split}
&2\rho\Vert \textbf{e}_n-\textbf{e}_{n-1}\Vert_{L^2}^2+\nu\Delta t\left(\Vert \nabla\textbf{e}_n\Vert_{L^2}^2-\Vert \nabla\textbf{e}_{n-1}\Vert_{L^2}^2+\Vert \nabla(\textbf{e}_n-\textbf{e}_{n-1})\Vert_{L^2}^2\right)\\
&\leq 2\Delta t|(r_n,\textbf{e}_n-\textbf{e}_{n-1})|+2\Delta t|(((\textbf{v}_{U_m})_n\cdot \nabla)\textbf{e}_n,\textbf{e}_n-\textbf{e}_{n-1})|+2\Delta t|((\textbf{e}_n\cdot \nabla)\textbf{v}_{U_m}(t_n),\textbf{e}_n-\textbf{e}_{n-1})|\\
&\leq \rho\Vert \textbf{e}_n-\textbf{e}_{n-1}\Vert_{L^2}^2+C(\Delta t)^2\Vert r_n\Vert_{L^2}^2+C(\Delta t)^2\Vert (\textbf{v}_{U_m})_n\Vert_{H^2}^2\Vert \nabla \textbf{e}_n\Vert_{L^2}^2+C(\Delta t)^2\Vert \nabla \textbf{e}_n\Vert_{L^2}^2,
\end{split}
\end{equation}
where we use the Sobolev embedding theorem $H^2\hookrightarrow L^{\infty}$.

Summing up (\ref{the3.11}) for $n$ from $1$ to $N$, we obtain
\begin{equation}\label{the3.12}
\begin{split}
&\Delta t\Vert \nabla\textbf{e}_n\Vert_{L^2}^2+\sum\limits_{n=1}^{N}\Vert \textbf{e}_n-\textbf{e}_{n-1} \Vert_{L^2}^2+\Delta t\sum\limits_{n=1}^{N}\Vert \nabla(\textbf{e}_n-\textbf{e}_{n-1})\Vert_{L^2}^2\\
&\leq C(\Delta t)^3+C(\Delta t)^2\sum\limits_{n=1}^{N}\Vert(\textbf{v}_{U_m})_n\Vert_{H^2}^2\Vert \nabla \textbf{e}_n\Vert_{L^2}^2+C(\Delta t)^2\sum\limits_{n=1}^{N}\Vert \nabla \textbf{e}_n\Vert_{L^2}^2.
\end{split}
\end{equation}

By using (\ref{the3.10}) and the discrete Gronwall inequality, we have
\begin{equation}\label{the3.13}
\begin{split}
&\Delta t\Vert \nabla\textbf{e}_n\Vert_{L^2}^2+\sum\limits_{n=1}^{N}\Vert \textbf{e}_n-\textbf{e}_{n-1} \Vert_{L^2}^2+\Delta t\sum\limits_{n=1}^{N}\Vert \nabla(\textbf{e}_n-\textbf{e}_{n-1})\Vert_{L^2}^2\leq C(\Delta t)^3.
\end{split}
\end{equation}

From (\ref{the3.4}) and the regularity of the solution to the Stokes problem, with the help of (\ref{the3.10}) and (\ref{the3.13}), we obtain
\begin{equation}\label{the3.14}
\begin{split}
&\Delta t\sum\limits_{n=1}^{N}\Vert \textbf{e}_n\Vert_{H^2}^2\leq C\Delta t\sum\limits_{n=1}^{N}\Vert r_n\Vert_{L^2}^2+C\Delta t\sum\limits_{n=1}^{N}\left\Vert\left(\frac{\textbf{e}_n-\textbf{e}_{n-1}}{\Delta t}\right)\right\Vert_{L^2}^2\\
&+C\Delta t\sum\limits_{n=1}^{N}\Vert(\textbf{v}_{U_m})_n\Vert_{H^2}^2\Vert \nabla \textbf{e}_n\Vert_{L^2}^2+C\Delta t\sum\limits_{n=1}^{N}\Vert \textbf{e}_n\Vert_{L^2}^2\leq C(\Delta t)^2.
\end{split}
\end{equation}
The proof is completed.
\end{proof}
\end{lemma}

\begin{lemma}
\label{lemma2}
(Lemma 20 in \cite{frei2020efficient} and following arguments of Lemma 3.4 in \cite{wang2022fast}). Let $0\leq u_1, u_2\leq u_{max}$, $\textbf{v}_1, \textbf{v}_2\in H^2(\Omega)$. The reaction term  (\ref{the1.2}) is Lipschitz continuous with respect to both arguments
\begin{equation}\label{the3.15}
\begin{split}	
|R(\textbf{v}_1,u_1)-R(\textbf{v}_2,u_2)|\leq C(\Vert \textbf{v}_1-\textbf{v}_2\Vert_{H^2}+|u_1-u_2|).
\end{split}
\end{equation}
Let the concentration $U$ be fixed, $\textbf{v}_U$ is the solution of the time periodic problem \ref{problem2}, it hols that 
\begin{equation}\label{the3.16}
\begin{split}	
\left\Vert \frac{\partial \textbf{v}_U}{\partial U}\right\Vert_{L^2}^2+\int_{0}^1\left\Vert \frac{\partial \textbf{v}_U}{\partial U}\right\Vert_{H^2}^2 dt\leq C.
\end{split}
\end{equation}
\end{lemma}

\begin{theorem}
\label{theorem}
Let $U\in C^2[0,T]$ and $U_m$ be the solutions to the auxiliary problem \ref{problem2} and the discrete equation (\ref{the2.2}), respectively. For $T_m=O(\epsilon^{-1})$, we obtain the following error estimate
\begin{equation}\label{the3.17}
\begin{split}	
|U_m-U(T_m)|=O(\Delta t+\epsilon\Delta T).
\end{split}
\end{equation}
\begin{proof}
Differentiating the equation (\ref{the2.1c}) with respect to $t$, and using Lemma \ref{lemma2}, we have
\begin{equation}\label{the3.18}
\begin{split}	
|U''(t)|&=\left|R\left(\textbf{v}_{U(t)}(1),U(t)\right)-R\left(\textbf{v}_{U(t)}(0),U(t)\right)+\epsilon\int_0^1 d_t R\left(\textbf{v}_{U(t)}(s),U(t)\right)\right|ds\\
&\leq C\epsilon\left(|U'(t)|+|U'(t)|\int_0^1 \left\Vert\frac{\partial \textbf{v}_U}{\partial U}\right\Vert_{H^2}ds\right)=O(\epsilon^2).
\end{split}
\end{equation}

Let $E_m=|U_m-U(T_m)|$, using Taylor expansion of $U(T_{m})$ around $T_{m-1}$ and (\ref{the2.2}), we obtain the error equation:
\begin{equation}\label{the3.19}
\begin{split}	
|E_m|&=\left|E_{m-1}+\Delta T\Delta t\epsilon\sum\limits_{n=1}^{N}\left(R((\textbf{v}_{U_{m-1}})_n,U_{m-1})-R(\textbf{v}_{U(T_{m-1})}(t_n),U(T_{m-1}))\right) \right|\\
&+\epsilon \Delta T O(\Delta t)+O((\Delta T)^2)|U''(\xi)|, \ \xi\in(T_{m-1},T_{m}).
\end{split}
\end{equation}

With the help of Lemma \ref{lemma2} and the Cauchy–Schwarz inequality, we have
\begin{equation}\label{the3.20}
\begin{split}	
&|R((\textbf{v}_{U_{m-1}})_n,U_{m-1})-R(\textbf{v}_{U(T_{m-1})}(t_n),U(T_{m-1}))|\\
&\leq C\left(\left\Vert(\textbf{v}_{U_{m-1}})_n-\textbf{v}_{U_{m-1}}(t_n)\right\Vert_{H^2}
+\left|\textbf{v}_{U_{m-1}}(t_n)-\textbf{v}_{U(T_{m-1})}(t_n)\right\Vert_{H^2}+|E_{m-1}|\right)\\
&\leq C\left(\left\Vert(\textbf{v}_{U_{m-1}})_n-\textbf{v}_{U_{m-1}}(t_n)\right\Vert_{H^2}
+|E_{m-1}|\right).
\end{split}
\end{equation}

By (\ref{the3.18})-(\ref{the3.20}) and using Lemma \ref{lemma}, we obtain
\begin{equation}\label{the3.21}
\begin{split}	
|E_m|&\leq |E_{m-1}|+C\epsilon\Delta T\left(|E_{m-1}|+\Delta t+\epsilon\Delta T\right)
\end{split}
\end{equation}
with a constant $C$ does not depend on $\epsilon$, $\Delta t$ and $\Delta T$. 

Summing over $m=1,\cdots,M$, with $E_{0}=0$, we obtain
\begin{equation}\label{the3.22}
\begin{split}	
|E_M|&\leq  C\epsilon\Delta T\sum\limits_{m=1}^{M-1}|E_{m}|+CT_{M}\epsilon(\Delta t+\epsilon\Delta T).
\end{split}
\end{equation}

Applying the discrete Gronwall inequality, for $T=O(\epsilon^{-1})$, we can easily obtain
\begin{equation}\label{the3.23}
\begin{split}	
|E_M|&\leq  C(\Delta t+\epsilon\Delta T)e^{C\epsilon T_M}\leq C(\Delta t+\epsilon\Delta T),
\end{split}
\end{equation}
which completes the proof.
\end{proof}
\end{theorem}

\section{Numerical experiment}
\label{section4}
In this section, we test a numerical example of 2-D Navier-Stokes problem, which is implemented based on the open-source finite element Software FreeFeM++. We choose the solution with fine macro and micro time steps as the reference to measure the error since we do not have possession of the exact solution. The shape function is assumed to be $\gamma(u,x)=ue^{-x^2}$, and the flow domain $\Omega=\left\{(x,y): |x|\leq 5, |y|\leq (2-ue^{-x^2})\right\}$. The spatial domain is discretized with 426 elements. The time periodic Dirichlet condition $\textbf{v}_{\text{in}}=20\left(1-\frac{y^2}{4}\right)\sin^2(\pi t)$ is set on the left inflow boundary, which makes the blood flow appear periodic. On the outflow boundary we set the pressure condition: $-p\vec{n}+\rho\nu\frac{\partial \textbf{v}}{\partial \vec{n}}=0$. In addition, we set computation time $T=4.8\cdot 10^4$, $u_0=0$, $\rho=1$, $\nu=0.04$ and $\sigma_0=30$. 

We test the accuracy of the proposed schemes (\ref{the2.2})-(\ref{the2.3}). We first consider the reference solutions which are obtained by using a fine micro time step $\Delta t^f=1/256$ to test the convergence rate of $\Delta t$. A set of decreasing micro time step $\Delta t=32\Delta t^{f}$, $\Delta t=16\Delta t^{f}$, $\Delta t=8\Delta t^{f}$, and $\Delta t=4\Delta t^{f}$ is used to perform the numerical simulation until $T=4.8\cdot 10^{5}$. Further, we test the convergence rate of $\Delta T$ by using the same method. The convergence rate of $\epsilon$ is tested by setting fine macro time steps. These numerical results are shown in Table \ref{table1}, which are consistent with the error estimates in Theorem \ref{theorem}. 

\begin{figure}[H]
	\centering
	\subfigure[$t=1.6\cdot10^{4}$]{
		\includegraphics[scale=0.24]{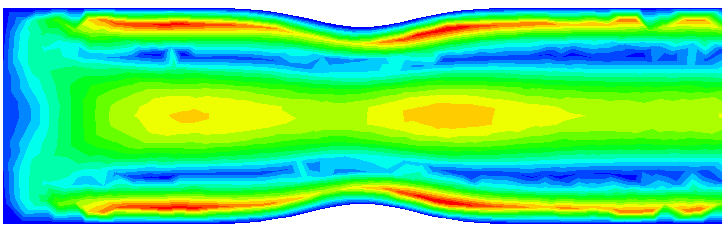}
	}
	\subfigure[$t=3.2\cdot10^{4}$]{
		\includegraphics[scale=0.24]{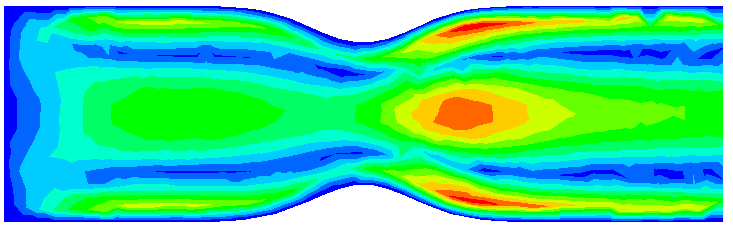}
		
	}
	\subfigure[$t=4.8\cdot10^{4}$]{
	\includegraphics[scale=0.24]{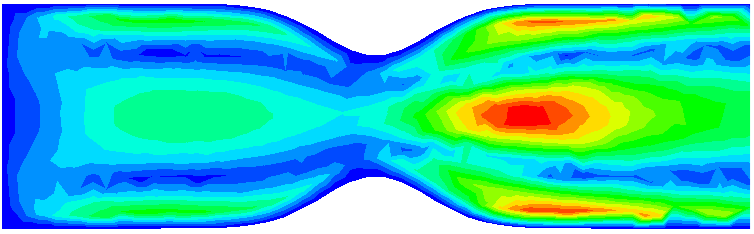}
	
}
	
	\caption{Snapshots of velocity field and plaque growth with $\Delta=1/32$, $\Delta T=2000$ and $\epsilon=2\cdot10^{-4}$.}
	\label{figure3}
\end{figure}

For the long-term simulation of this plaque growth problem, we perform direct computation by using a time step of $\Delta t=1/20$ that takes 71 hours. It is easy to see from Table \ref{table1} that the effective numerical multiscale method we used significantly reduces the computational time, and it increases as the macro and micro steps become finer. Figure \ref{figure3} shows that the plaque grows significantly and a strong narrowing of the flow domain with the increase in time.

\begin{table}[H]
	{\footnotesize
		\caption{Errors of macro variable $U$ at $t=T=4.8\cdot 10^{4}$, convergence rates of $\Delta t$, $\Delta T$ with fixed $\epsilon=2\cdot 10^{-4}$, and convergence rates of $\epsilon$.}\label{table1}
		\resizebox{\linewidth}{!}{
			\begin{tabular}{c c c c} \hline  		
				$\Delta t$ ($\Delta T=500$)   & Error  & Order & CPU time(s) \\ \hline
				1/8		& 5.181e-3 &  & 2318.0   \\
				1/16	& 2.645e-3 & 0.97 & 2870.4   \\  
				1/32	& 1.251e-3  &  1.08 & 3042.0  \\  
				1/64	& 6.169e-4  & 1.02 & 3311.6   \\  \hline
				$\Delta T$ ($\Delta t=1/32$)	&   &  &    \\  \hline
				8000		& 2.403e-2 &  & 104.3   \\
				4000	& 1.094e-2    & 1.14 & 242.2   \\  
				2000	&  5.082e-3   &  1.11 & 563.9  \\  
				1000	& 2.501e-3  &   1.02 & 1124.2  \\  \hline
				$\epsilon$ ($\Delta T=2000$,$\Delta t=1/32$)	&   &  &    \\  \hline
				$4\cdot10^{-4}$		&  9.750e-3 &  &  555.2  \\
				$2\cdot10^{-4}$	& 5.082e-3     & 0.94 & 563.9    \\  
				$1\cdot10^{-4}$	&  2.594e-3   &  0.97 &  563.2 \\  
				$5\cdot10^{-5}$	& 1.297e-3  &   1.00 &  588.2 \\  \hline
		\end{tabular}}
	}
\end{table}

\section*{Acknowledgments}
This work is supported by the National Natural Science Foundation of China (Nos. 11861131004, 11771040).

\bibliographystyle{elsarticle-num}

\bibliography{reference}

\end{document}